\documentclass[11pt,twoside,reqno]{amsart}
\usepackage{amsmath,times,amssymb,amsbsy,amscd,amsfonts,amstext,graphicx}
\calclayout 
\newtheorem{thm}{\bf Theorem}
\newtheorem{lem}{\bf Lemma}
\newtheorem{prop}{\bf Proposition}
\newtheorem{cor}{\bf Corollary}
\newtheorem{rk}{\bf Remark}

\def\n{\noindent}

\begin{document}
\title
{UNIVERSAL UPPER BOUND FOR THE  GROWTH OF ARTIN MONOIDS}
\author{Barbu Berceanu$^{*,**}$, Zaffar Iqbal$^{*}$}
\thanks{{ $^{*}$ Abdus Salam School of
Mathematical Sciences, GC University, Lahore, Pakistan, Email:
zaffarsms@yahoo.com}
{\indent \!$^{**}$} Simion Stoilow Institute of Mathematics, Bucharest,
Romania, Email: Barbu.Berceanu@imar.ro}
 \fontsize{11}{21.3}\selectfont
\begin{abstract}
In this paper we study the growth rates of Artin monoids and we show
that $4$ is a universal upper bound. We also show that the
generating functions of the associated right-angled Artin monoids
are given by families of Chebyshev
polynomials. Applications to Artin groups and positive braids are given.
\end{abstract}
\maketitle
\pagestyle{myheadings} \markboth{\centerline {\scriptsize B. Berceanu, Z.
Iqbal}} {\centerline {\scriptsize Universal bound for the growth of
Artin monoids}}

\section*{\bf \S\,1. Introduction}

One can prove that $b_{k}=$ the number of  braids with three strands
and $k$ positive crossings is given by $b_{0}=1$, $b_{1}=2$,
$b_{2}=4$, $b_{3}=7$, $b_{4}=12$, $b_{5}=20$; in general we have
$b_{k}=F_{k+2}-1,$ where $\{F_{k}\}$ are Fibonacci numbers.
Therefore the growth function of this sequence is exponential with
base\, $\frac{1+\sqrt{5}}{2}$ .

In this paper we generalize this to all Artin spherical monoids (and
also to some associated right-angled Artin monoids) and we find a
universal upper bound for the growth functions of all these monoids.

An \emph{Artin spherical monoid (or group)} is given by a finite
union of connected Coxeter graphs from the classical list of Coxeter
diagrams (see \cite{bbk}; also, for a recent survey, see
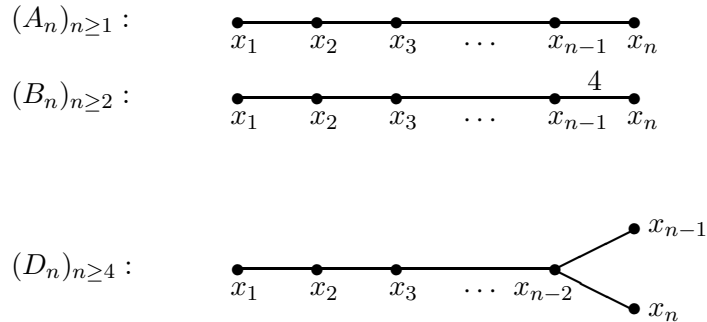
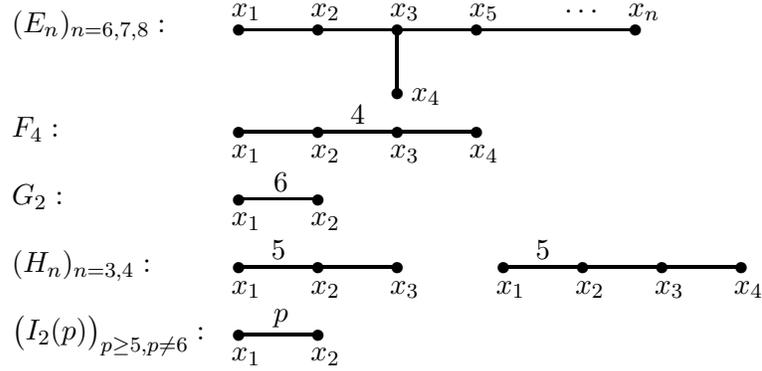
\cite{paris}):

\begin{center}
\begin{picture}(300,40)
\thicklines \put(100,27){\line(1,0){150}}
\put(98,24){$\bullet$}\put(128,24){$\bullet$}\put(158,24){$\bullet$}
\put(218,24){$\bullet$}\put(248,24){$\bullet$}
\put(98,17){$x_{1}$}\put(128,17){$x_{2}$}\put(158,17){$x_{3}$}
\put(218,17){$x_{n-1}$}\put(248,17){$x_{n}$}\put(186,16){$\cdots$}
\put(15,25){$(A_{n})_{n\geq1}:$}
\end{picture}
\end{center}

\begin{center}
\begin{picture}(300,24)
\thicklines \put(100,27){\line(1,0){150}}
\put(98,24){$\bullet$}\put(128,24){$\bullet$}\put(158,24){$\bullet$}
\put(218,24){$\bullet$}\put(248,24){$\bullet$}
\put(98,17){$x_{1}$}\put(128,17){$x_{2}$}\put(158,17){$x_{3}$}
\put(218,17){$x_{n-1}$}\put(248,17){$x_{n}$}\put(186,16){$\cdots$}
\put(15,25){$(B_{n})_{n\geq2}:$} \put(233,30){$4$}
\end{picture}
\end{center}

\begin{center}
\begin{picture}(300,60)
\thicklines
\put(100,27){\line(1,0){120}}\put(220,26.5){\line(2,1){30}}\put(220.5,26.5){\line(2,-1){30}}
\put(98,24){$\bullet$}\put(128,24){$\bullet$}\put(158,24){$\bullet$}
\put(218,24){$\bullet$}\put(248,39.5){$\bullet$}\put(248,9){$\bullet$}
\put(98,17){$x_{1}$}\put(128,17){$x_{2}$}\put(158,17){$x_{3}$}
\put(205,17){$x_{n-2}$}\put(256,9){$x_{n}$}\put(256,41){$x_{n-1}$}\put(186,16){$\cdots$}
\put(15,25){$(D_{n})_{n\geq4}:$}
\end{picture}
\end{center}

\begin{center}
\begin{picture}(300,35)
\thicklines
\put(100,27){\line(1,0){150}}\put(160.5,27){\line(0,-1){24}}
\put(98,24){$\bullet$}\put(128,24){$\bullet$}\put(158,24){$\bullet$}
\put(188,24){$\bullet$}\put(248,24){$\bullet$}\put(158,0){$\bullet$}
\put(98,32){$x_{1}$}\put(128,32){$x_{2}$}\put(158,32){$x_{3}$}\put(166,0){$x_{4}$}
\put(188,32){$x_{5}$}\put(248,32){$x_{n}$}\put(224,31){$\cdots$}
\put(15,26){$(E_{n})_{n=6,7,8}:$}
\end{picture}
\end{center}

\begin{center}
\begin{picture}(300,35)
\thicklines \put(100,27){\line(1,0){90}}
\put(98,24){$\bullet$}\put(128,24){$\bullet$}\put(158,24){$\bullet$}
\put(188,24){$\bullet$}
\put(98,17){$x_{1}$}\put(128,17){$x_{2}$}\put(158,17){$x_{3}$}
\put(188,17){$x_{4}$}\put(143,30){$4$} \put(15,25){$F_{4}:$}
\end{picture}
\end{center}

\begin{center}
\begin{picture}(300,22)
\thicklines \put(100,27){\line(1,0){30}}
\put(98,24){$\bullet$}\put(128,24){$\bullet$}
\put(98,17){$x_{1}$}\put(128,17){$x_{2}$} \put(114,30){$6$}
\put(15,25){$G_{2}:$}
\end{picture}
\end{center}

\begin{center}
\begin{picture}(300,22)
\thicklines \put(100,27){\line(1,0){60}}
\put(98,24){$\bullet$}\put(128,24){$\bullet$}\put(158,24){$\bullet$}
\put(98,17){$x_{1}$}\put(128,17){$x_{2}$}\put(158,17){$x_{3}$}
\put(113,30){$5$}
\put(200,27){\line(1,0){90}}
\put(198,24){$\bullet$}\put(228,24){$\bullet$}\put(258,24){$\bullet$}\put(288,24){$\bullet$}
\put(198,17){$x_{1}$}\put(228,17){$x_{2}$}\put(258,17){$x_{3}$}\put(288,17){$x_{4}$}
\put(213,30){$5$}\put(15,25){$(H_{n})_{n=3,4}:$}
\end{picture}
\end{center}

\begin{center}
\begin{picture}(300,22)
\thicklines \put(100,27){\line(1,0){30}}
\put(98,24){$\bullet$}\put(128,24){$\bullet$}
\put(98,17){$x_{1}$}\put(128,17){$x_{2}$} \put(114,32){$p$}
\put(15,25){$\big(I_{2}(p)\big)_{p\geq5,p\neq6}:$}
\end{picture}
\end{center}

By convention $r_{ij}$ is the label of the
edge between $x_{i}$ and $x_{j}$ $(i\neq j);$ if there is no label then
$r_{ij}=3$. If there is no edge between $x_{i}$ and $x_{j},$ then
$r_{ij}=2$.

To a given diagram we associate the monoid $\mathcal{X}^{+}_{n}$ with the
following presentation (generators corresponds to the vertices, and
relations corresponds to the labels $r_{ij}$ of the graphs):
$$\mathcal{X}^{+}_{n}=\Big\langle
x_{1},x_{2},\ldots,x_{n}\,\big|\,\mathop{\underbrace{x_{i}x_{j}x_{i}x_{j}\cdots}}
\limits_{r_{ij}\,\textrm{times}}=\mathop{\underbrace{x_{j}x_{i}x_{j}x_{i}\cdots}}
\limits_{r_{ij}\,\textrm{times}}\Big\rangle;$$
the groups $\mathcal{X}_{n}$ are defined by the same presentation.

Again by convention, if $r_{ij}=\infty$ there is no relation between
$x_{i}$ and $x_{j}$. We obtain the following list of Artin monoids
(and groups $\mathcal{X}_{*}$) corresponding to connected graphs:
$$\mathcal{X}^{+}_{*}=\{A^{+}_{n}(n\geq1),B^{+}_{n}(n\geq2), D^{+}_{n}(n\geq4),
E^{+}_{n}(n=6,7,8),F^{+}_{4},G^{+}_{2},H^{+}_{3},H^{+}_{4},I^{+}_{2}(p)\}.$$
For example, the associated group (monoid)
to the graph $A_{n}$ is the classical \emph{Artin group (and monoid $A^{+}_{n}$) of
(positive) braids} (see \cite{art} and \cite{gar}):
$$A_{n}=\left\langle x_{1},x_{2},\cdots\!,x_{n}\,\Bigg|
   \begin{array}{l}
               \!\!x_{i}x_{j}=x_{j}\,x_{i}\,\,\,\mbox{if} \, \mid i-j\mid\,\,\geq2\\
               \!\!x_{i+1}\,x_{i}\,x_{i+1}=x_{i}\,x_{i+1}\,x_{i}\,\,\mbox
               {if}\,\, 1\leq i\leq n-1
   \end{array}\!\!\right\rangle.$$

We study the growth function of these monoids (or groups), with a
fixed set of generators $\{x_{i}\}_{i=1,\cdots,n}$ , i.e., the
growth of the sequence $k\mapsto \#$\{words of length $k$ in the
given alphabet\}.

Let us look at two extremal Coxeter graphs, say $N_{n}$ and $C_{n}$,
where $N_{n}$ is the totally disconnected graph with $n$ vertices
(all the labels are 2) and $C_{n}$ is the complete graph and all
labels are $\infty$. The first monoid $N^{+}_{n}$ is the free abelian monoid
with Hilbert series
$$\frac{1}{(1-t)^{n}}=\sum\limits^\infty_{k=0}
\small\begin{pmatrix}n-1+k\\n-1\\\end{pmatrix}t^{k}$$ with polynomial
growth function. The second monoid $C^{+}_{n}$  is the free monoid with Hilbert
series
$$\frac{1}{1-nt}=\sum\limits^\infty_{k=0} n^{k}t^{k}$$
and the growth function is exponential $n^{k}$  and this is
unbounded for the sequence $(C_{n})_{n\geq1}$.

In spite of the last example, we shall show that the family of
Artin spherical monoids has a universal upper bound for the growth function.

More precisely, for a sequence $\{s_{k}\}_{k\geq1}$ of positive numbers, we define the growth rate
by:

\noindent\textbf{Definition 1.} We say that $\{s_{k}\}_{k\geq1}$ has
a \emph{growth rate less than} $\gamma$ ($\gamma$ is a real number)
if
$$\overline{\lim\limits_k}\exp\Big(\frac{\log s_{k}}{k}\Big)<\gamma.$$
This is equivalent to the fact that there exists
$\gamma_{0}<\gamma$ and $c>0$ such that  $s_{k}<c\,\gamma_{0}^{k}$
for all $k$. For instance, if $\max\limits_{1\leq i\leq n}|\alpha_{i}|<\gamma$,
the sequence $c_{k}=\sum\limits^n_{i=1}\beta_{i} \alpha_{i}^{k}$ has
a growth rate less than $\gamma$. In this case we have also
$\frac{c_{k+1}}{c_{k}}<\gamma$ for $k\gg0$; if the sequence
$\big\{\frac{c_{k+1}}{c_{k}}\big\}_{k}$ is convergent, its limit is
called the \emph{growth ratio} of the sequence $\{c_{k}\}$. See
\cite{ufn} for a general discussion of growth functions.

In \cite{p.del} P. Deligne proved that all Artin spherical monoids
have rational Hilbert series and he gave a formula in terms of the
numbers of reflections in the associated Weyl groups.

In \cite{xu} P. Xu studied the growth function of the monoid of
positive braids (the series $A_{n}$) and she found explicit formulae
for the Hilbert series of $A_{3}$ and $A_{4}$.

The main result of this paper is
\begin{thm}\label{th1}
The growth rate of all Artin spherical monoids (with classical generators) is less than $4$.
\end{thm}
In order to prove these results we introduce the associated
\emph{right-angled} Artin spherical monoids:
$$\mathcal{X}^{\infty}_{*}=\{A^{\infty}_{*},B^{\infty}_{*},
D^{\infty}_{*},E^\infty_{*},F^{\infty}_{4},G^{\infty}_{2},
H^{\infty}_{3},H^{\infty}_{4},I^{\infty}_{2}(p)\}$$ where all the
labels $\geq3$ are replaced by $\infty$. For example, to the Artin
monoid of positive braids we associate the monoid $A^\infty_{n}$
given by the Coxeter graph:
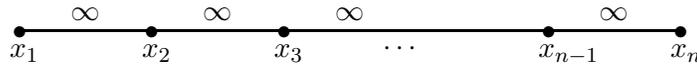
\begin{center}
 \begin{picture}(350,25)
 \thicklines
 \put(52,10){\line(1,0){250}}
 \put(50,7){$\bullet$} \put(100,7){$\bullet$}
 \put(150,7){$\bullet$} \put(250,7){$\bullet$} \put(300,7){$\bullet$}
 \put(49,0){$x_{1}$} \put(99,0){$x_{2}$}\put(149,0){$x_{3}$}\put(249,0){$x_{n-1}$}\put(299,0){$x_{n}$}
 \put(72,14){$\infty$}\put(122,14){$\infty$}\put(172,14){$\infty$}\put(272,14){$\infty$}\put(190,1){$\cdots$}
 \end{picture}
 \end{center}

The sequence $A^{\infty}_{*}$ contains $B^{\infty}_{*},$
$F^{\infty}_{4},$ $G^{\infty}_{2},$ $H^{\infty}_{*},$
$I^{\infty}_{2}(p),$ so the study of growth rate is reduced to two
sequences of monoids $(A^{\infty}_{*}$ and $D^{\infty}_{*})$ and
three exceptional cases $E^{\infty}_{6},E^{\infty}_{7},$ and
$E^{\infty}_{8}$. Next we reduce the series $A^{\infty}_{*}$ and
$D^{\infty}_{*}$ to a unique series: we introduce a new monoid
$K^{\infty}_{n}$ with Coxeter graph:
 
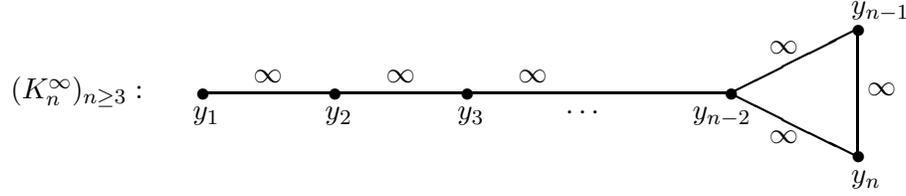
\begin{center}
 \begin{picture}(300,70)
 \thicklines
 \put(52,37){\line(1,0){200}}\put(252,37){\line(2,1){48}}
 \put(252,37){\line(2,-1){48}}\put(300.5,12){\line(0,1){48}}
 \put(50,34){$\bullet$} \put(100,34){$\bullet$}
 \put(150,34){$\bullet$} \put(250,34){$\bullet$} \put(298,58){$\bullet$}\put(298,10){$\bullet$}
 \put(49,27){$y_{1}$} \put(99,27){$y_{2}$}\put(149,27){$y_{3}$}\put(238,27){$y_{n-2}$}
 \put(298,67){$y_{n-1}$}\put(298,3){$y_{n}$}
 \put(-20,35){$(K^{\infty}_{n})_{n\geq3}:$}
 \put(72,41){$\infty$}\put(122,41){$\infty$}\put(172,41){$\infty$}
 \put(267,52){$\infty$}\put(267,18){$\infty$}\put(304,36){$\infty$}
 \put(190,27){$\cdots$}
 \end{picture}
 \end{center}
and presentation
$$K^{\infty}_{n}=\left\langle
y_{1},y_{2},\ldots,y_{n}\,\Bigg|
   \begin{array}{l}
               y_{i}y_{j}=y_{j}\,y_{i}\,,\,\,\,\mbox{} \, j+2\leq i \leq\,n-1\\
               y_{n}y_{k}=y_{k}y_{n}\,,\,\,\mbox
               {}\,\, 1\leq k \leq\,n-3
   \end{array}\!\right\rangle.$$
$\big($Because we will not use the groups associated to Coxeter graphs with labels $\infty$,
we simplify the notations:
$A^{\infty}_{n},\cdots,K^{\infty}_{n}$ instead of $(A^{\infty}_{n})^{+},\cdots,(K^{\infty}_{n})^{+}.\big)$

The growth rate of $\mathcal{X}^{+}_{*}$ ($\mathcal{X}_{*}$ in the classical list, but not
$E_{*}$) is less than the growth rate of
$K^{\infty}_{*}$ as a consequence of the next proposition:
\begin{prop}\label{p1}
There exist canonical surjective homogeneous maps $:$
$$K^{\infty}_{n}\mathop{\twoheadrightarrow}\limits^\phi\mathcal{X}^{\infty}_{n}
\mathop{\twoheadrightarrow}\limits^\psi\mathcal{X}^{+}_{n}.$$
\end{prop}

Growth functions for Artin groups were studied by many authors, using
different set of generators. J. Mairesse  and  F. Math\'{e}us \cite{mair} proved that
the dihedral series $I_{2}(p)$ with classical generators
$\{x_{1},x_{2},x^{-1}_{1},x^{-1}_{2}\}$ has a rational growth function;
one can easily show that $4$ is still an upper bound for the
growth rate. R. Charney \cite{char} analyzed the growth function of Artin groups using
normal generators $\{d,d^{-1}: \mbox{d is a divisor of the
Garside element}\,\,\Delta\}$ and the same groups $\{I_{2}(p)\}$ have
unbounded growth ratio $p-1$. In a recent paper \cite{mang} J. Mangahas
shows that the growth rates of braid groups have a lower bound which does not depend on
the (finite) set of generators.

The main idea was to reduce the study of Artin groups $\mathcal{X}_{n}$ to the study of the associated
monoid $\mathcal{X}^{+}_{n}$. F. A. Garside \cite{gar} and P. Deligne \cite{p.del}
used the fundamental word $\Delta_{n}\in\mathcal{X}^{+}_{n}$ and the canonical decomposition
$$\mathcal{X}_{n}=\bigcup\limits_{k\in\mathbb{Z}}\Delta_{n}^{k}\,\mathcal{X}^{+}_{n};$$
using these one can show
\begin{thm}\label{th2}
All Artin spherical groups $\mathcal{X}_{n}$, as monoids generated by
$\{x_{1},\cdots,x_{n},\nabla=\Delta_{n}^{-1}\}$
have growth rate less than $4$.
\end{thm}

In $\S \,2$ we solve the word problem for $K^{\infty}_{n}$ and we compute
its Hilbert series and also those of
$A^{\infty}_{n}$ and $D^{\infty}_{n}$.

In $\S\,3$ we show that $K^{\infty}_{n}$ has growth rate less than 4
using a sequence of Chebyshev type polynomials.
We give in proposition \ref{Rn} sufficient conditions for Chebyshev type sequences of
polynomials to have only real roots in a bounded interval.

In $\S\,4$ we describe the results corresponding to the exceptional
monoids $E^{\infty}_{6},E^{\infty}_{7}$ and $E^{\infty}_{8}$.

The last paragraph contains the proofs of the theorems \ref{th1} and \ref{th2} and of the
proposition \ref{p1} and also the extension of all the previous
results to the Artin spherical monoids with non-connected Coxeter
graphs.

As an application to braids, we have the following
\begin{cor}
If $b^{[n]}_{k}$ is the number of positive $n$-braids with $k$
positive crossings, then for any $n$ and for large values of $k$ we
have $b^{[n]}_{k+1}<4\,b^{[n]}_{k}$.
\end{cor}

\section*{\bf \S \,2.\,The Hilbert series of the monoid $K^{\infty}_{n}$}

In a presentation of a monoid we fix a total order of the
generators; in all our examples we choose the natural order
$y_{1}<y_{2}<\cdots<y_{n}$. Such a presentation is \emph{complete}
if and only if all the ambiguities are solvable (see \cite{berg},
\cite{coh}).
\begin{lem}\label{cpc}
 The presentation
$$K^{\infty}_{n}=\left\langle y_{1},y_{2},\ldots,y_{n}\,\Bigg|
   \begin{array}{l}
               y_{i}y_{j}=y_{j}\,y_{i}\,,\,\,\,\mbox{} \, j+2\leq i \leq\,n-1\\
               y_{n}y_{k}=y_{k}y_{n}\,,\,\,\mbox
               {}\,\, 1\leq k \leq\,n-3
   \end{array}\!\right\rangle$$ is a complete presentation.
\end{lem}

\proof All the  ambiguities are of type
$y_{i}\,y_{j}\,y_{k}\,,\,i-2\geq j\geq k+2$ ($i=n$ should be $\geq
j+3$) and these can be simplified in two ways; after three steps we
get the canonical word $y_{k}\,y_{j}\,y_{i}$.\endproof

As a consequence we obtain a solution of the word problem in
$K^{\infty}_{n}$. The smallest words in length-lexicographic order
are given by the next proposition and the defining relations are
sufficient for the rewriting system of this monoid:

\begin{prop}\label{can.forms}
In $K^{\infty}_{n}$ the canonical form of a non empty word is given
by
$$y^{a_{0}}_{i_{1}}y^{a_{1}}_{i_{1}-1}\cdots
y^{a_{k_{1}}}_{i_{1}-k_{1}}y^{b_{0}}_{i_{2}}y^{b_{1}}_{i_{2}-1}\cdots
y^{b_{k_{2}}}_{i_{2}-k_{2}}\cdots
y^{h_{0}}_{i_{p}}y^{h_{1}}_{i_{p}-1}\cdots
y^{h_{k_{p}}}_{i_{p}-k_{p}}$$ where $i_{j}-k_{j}<i_{j+1}$,
$j=1,\ldots,p-1$ and all the exponents are positive; as a special case,
 after the factor $y^{h_{*}}_{n}$ could also appear $y^{j_{*}}_{n-2}$ .
\end{prop}

\begin{rk}
Similar results are true for $A^{\infty}_{*},D^{\infty}_{*},$ and
$E^{\infty}_{*}$.
\end{rk}
Now we start to compute the Hilbert series
$\mathcal{\mathcal{H}}^{[n]}_K(t)=\sum\limits_{k\geq0}
c^{[n]}_{k}t^{k}$ of $K^{\infty}_{n}$ and also
$\mathcal{H}^{[n]}_{K;i}(t)=\sum\limits_{k\geq0}
c^{[n]}_{k;i}t^{k}$; here $c^{[n]}_{k}$ (or $c_{k}$) and
$c^{[n]}_{k;i}$ (or $c_{k;i}$) denote the number of words in
$K^{\infty}_{n}$ of length $k,$ respectively words (in canonical
form) of length $k$ starting with $y_{i}$. In the same way we denote
by $\sum\limits_{k\geq0} a_{k}t^{k},$ $\sum\limits_{k\geq0}
b_{k}t^{k}$ the Hilbert series (or generating functions) of the
monoids $\mathcal{X}^{+}_{*}$, respectively
$\mathcal{X}^{\infty}_{*}$.

\begin{cor}\label{Ksys.rr}
In the monoid $K^{\infty}_{n}$ the following relations are satisfied

$a)$ $c_{0}=1,$ $c_{1;i}=1,$ $c_{k} =\sum\limits_{i=1}^n
c_{k;i}$\,\, $(k\geq1)$ .

$b)$ $c_{k;i}\,(k\geq2)$ are given by the recurrence
 $$\left\{
  \begin{array}{ll}
    c_{k;j}=\sum\limits_{i=j-1}^n c_{k-1;i}\,\,(j\neq n),\\
    c_{k;n}=\sum\limits_{i=n-2}^n c_{k-1;i}\,.
  \end{array}
\right.$$
\end{cor}

The characteristic polynomial of this recurrence is given by:
$$\mathcal{K}_{n}(\lambda)= \det\begin{bmatrix}
  \lambda-1 & -1 & \cdots &  -1 & -1 & -1 &  -1 \\
  -1 & \lambda-1 & \cdots& -1 & -1 & -1 &  -1 \\
   0 & -1 & \cdots & -1 &  -1 & -1 & -1\\
  \vdots  & \vdots && \vdots  & \vdots & \vdots & \vdots \\
   0 & 0  & \cdots & -1 & \lambda-1& -1  &  -1 \\
   0 & 0  & \cdots & 0  &  -1 & \lambda-1 & -1 \\
   0 & 0  & \cdots & 0  &  -1 & -1 & \lambda-1 \\
\end{bmatrix}$$
\begin{lem}
The polynomials $\big(\mathcal{K}_{n}(\lambda)\big)_{n\geq3}$
satisfy the recurrence$:$
\begin{equation}\label{Snrr}
    \mathcal{K}_{n}(\lambda)=\lambda \mathcal{K}_{n-1}(\lambda)-\lambda
    \mathcal{K}_{n-2}(\lambda)\,\,\,\,\,(n\geq5)
\end{equation}
with $\mathcal{K}_{3}(\lambda)=\lambda^3-3\lambda^2,
\mathcal{K}_{4}(\lambda)=\lambda^4-4\lambda^3+2\lambda^2$.
\end{lem}
\proof Decompose $\mathcal{K}_{n}$ as a sum of two determinants
$U_{n}$ and $V_{n}$ with the first rows given by
$[\lambda,0,\ldots,0]$ and $[-1,-1,\ldots,-1]$ respectively.
Elementary operations on the first two rows will give
$U_{n}(\lambda)=\lambda \mathcal{K}_{n-1}(\lambda)$ and
$V_{n}(\lambda)=-\lambda \mathcal{K}_{n-2}(\lambda)$. We can extend
the recurrence for $n\geq2$ by defining
$\mathcal{K}_{0}(\lambda)=2,\mathcal{K}_{1}(\lambda)=\lambda,$ and
$\mathcal{K}_{2}(\lambda)=\lambda^2-2\lambda$.
\endproof
From corollary \ref{Ksys.rr} we have:

\begin{cor}The Hilbert series of $K^{\infty}_{n}$ is given by$:$

$$\left\{
  \begin{array}{ll}
    \mathcal{H}^{[n]}_{K}(t)\,\,= 1+\sum\limits_{i=1}^n \mathcal{H}^{[n]}_{K;i}(t)\,\,;  \\
    \mathcal{H}^{[n]}_{K;j}(t)= t+t\sum\limits_{i=j-1}^n \mathcal{H}^{[n]}_{K;i}(t)\,\,\,\,\,(j\neq1,n); \\
    \mathcal{H}^{[n]}_{K;j}(t)=\mathcal{H}^{[n]}_{K;j+1}(t)\,\,\,\,\,(j=1,\,n-1).
  \end{array}
\right.$$
\end{cor}
The system given by the last $n$ equations has the determinant:
$$\det W_{n}=t^{n}\mathcal{K}_{n}\Big(\frac{1}{t}\Big).$$

The first two algebras have Hilbert series:

$\mathcal{H}^{[3]}_{K}(t)=\frac{1}{1-3t}=1+3t+9t^2+27t^3+81t^4+\cdots$
\,\,and

$\mathcal{H}^{[4]}_{K}(t)=\frac{1}{1-4t+2t^2}=1+4t+14t^2+48t^3+164t^4+560t^5+\cdots.$

The characteristic polynomials $\mathcal{A}_{n}(\lambda)$ of
$A^{\infty}_{n}$ satisfy the same recurrence:
\begin{lem}
\begin{equation}\label{Pnrr}
    \mathcal{A}_{n}(\lambda)=\lambda \mathcal{A}_{n-1}(\lambda)-\lambda
    \mathcal{A}_{n-2}(\lambda)\,\,\,\,\,(n\geq2)
\end{equation}
with $\mathcal{A}_{0}(\lambda)=1$ and $\mathcal{A}_{1}(\lambda)=\lambda-1$.
\end{lem}
and also
\begin{lem}
The polynomials $\mathcal{K}_{n}(\lambda)$ and
$\mathcal{A}_{n}(\lambda)$ $(n\geq0)$ satisfy
\begin{equation}\label{SPrr}
\mathcal{K}_{n}(\lambda)=\lambda \mathcal{A}_{n-1}(\lambda)-\lambda^2 \mathcal{A}_{n-3}(\lambda)\,\,\,\,(n\geq3).
\end{equation}
\end{lem}
\proof This relation is true for $n=3$ and $n=4$ and
$\mathcal{A}_{n}$ and $\mathcal{K}_{n}$ satisfy the same
recurrence.\endproof

\begin{lem}\label{Hm}
The Hilbert series $\mathcal{H}^{[n]}_{K;m}(t),$ $2\leq m\leq n-1$
is given by
$$\mathcal{H}^{[n]}_{K;m}(t)=\frac{t^{m-1}\mathcal{A}_{m-2}(\frac{1}{t})}{\det W_{n}}\,.$$
\end{lem}
\proof  The $\det L_{m}$ in the Cramer's rule can be reduced to a
product of two determinants, one equal to
$t^{m-1}\mathcal{A}_{m-2}(\frac{1}{t})$ and a second one, given by a
unipotent $(n-m)$ triangular matrix.
\endproof

\begin{prop}\label{HtnSn} The Hilbert series of the monoid $K^{\infty}_{n}$
is given by  $$\mathcal{H}^{[n]}_{K}(t)=\frac{1}{t^n
\mathcal{K}_{n}(\frac{1}{t})}\,\,.$$
\end{prop}
\proof From the lemma \ref{Hm} we have
\begin{eqnarray*}
\mathcal{H}^{[n]}_{K}(t) &=&
1+\mathcal{H}^{[n]}_{K;1}(t)+\mathcal{H}^{[n]}_{K;2}(t)
+\cdots+\mathcal{H}^{[n]}_{K;n}(t)\\&=&
\frac{1}{t^n\mathcal{K}_{n}(\frac{1}{t})}\Big[t^{n}\mathcal{K}_{n}\Big(\frac{1}{t}\Big)+2t+
t^2\mathcal{A}_{1}\Big(\frac{1}{t}\Big)+t^3\mathcal{A}_{2}\Big(\frac{1}{t}\Big)+\cdots+
2t^{n-2}\mathcal{A}_{n-3}\Big(\frac{1}{t}\Big)\Big]\\&=&
\frac{1}{t^n\mathcal{K}_{n}(\frac{1}{t})}\Big[2t+
t^2\mathcal{A}_{1}\Big(\frac{1}{t}\Big)+t^3\mathcal{A}_{2}\Big(\frac{1}{t}\Big)+\cdots+
t^{n-2}\mathcal{A}_{n-3}\Big(\frac{1}{t}\Big)+t^{n-1}\mathcal{A}_{n-1}\Big(\frac{1}{t}\Big)\Big]\\
&=&\frac{1}{t^n\mathcal{K}_{n}(\frac{1}{t})}\Big[2t+t^2\mathcal{A}_{2}\Big(\frac{1}{t}\Big)\Big]=
\frac{1}{t^n\mathcal{K}_{n}(\frac{1}{t})} \indent\indent(\mbox{using
the relations \ref{Pnrr} and \ref{SPrr}})\,.
\end{eqnarray*}\endproof

Now we will separate the zero roots of $\mathcal{K}_{n}$ from the
others:
\begin{prop}\label{SnTn}
The polynomial $\mathcal{K}_{n}(\lambda)$ has the following form:
$$\mathcal{K}_{n}(\lambda)=\lambda^{[\frac{n+1}{2}]} K_{n}(\lambda)$$ where
$K_{n}$ is a polynomial of degree $[\frac{n}{2}];$ the sequence
$(K_{n})_{n\geq0}$  is defined by

$1)$ $K_{0}=2,K_{1}=1, K_{2}=\lambda-2, K_{3}=\lambda-3\,;$

$2)$ $K_{2p+1}=(\lambda-2)K_{2p-1}-K_{2p-3}$ and
$K_{2p+2}=(\lambda-2)K_{2p}-K_{2p-2}.$
\end{prop}

\proof We prove these relations by induction:
$$\mathcal{K}_{2p}=\lambda( \mathcal{K}_{2p-1}- \mathcal{K}_{2p-2})=
\lambda^{p}(\lambda K_{2p-1}-K_{2p-2})=\lambda^{p}K_{2p}$$ and
similarly
$$\mathcal{K}_{2p+1}=\lambda( \mathcal{K}_{2p}- \mathcal{K}_{2p-1})=
\lambda^{p+1}(\lambda K_{2p}-K_{2p-1})=\lambda^{p+1}K_{2p+1}.$$

We will check only the recurrence for odd polynomials:
\begin{eqnarray*}
K_{2p+1}&=& K_{2p}-K_{2p-1}=(\lambda K_{2p-1}-K_{2p-2})-K_{2p-1} \\
&=&(\lambda-1)K_{2p-1}-K_{2p-1}-K_{2p-3}=(\lambda-2)K_{2p-1}-K_{2p-3}.
\end{eqnarray*}
\endproof

\begin{rk}\label{rk2}
 The results of proposition \ref{HtnSn} hold for the monoids
$A^{\infty}_{n}$ and $D^{\infty}_{n}:$

$\textbf{A})$ The polynomial $\mathcal{A}_{n}(\lambda)$ has the form
$\mathcal{A}_{n}(\lambda)=\lambda^{[\frac{n}{2}]} A_{n}(\lambda)$
where $A_{n}$ is a polynomial of degree $[\frac{n+1}{2}];$ the
sequence $(A_{n})_{n\geq-1}$ is defined by:

$1)$ $A_{-1}=1,A_{0}=1, A_{1}=\lambda-1,A_{2}=\lambda-2.$

$2)$ $A_{2p+1}=(\lambda-2)A_{2p-1}-A_{2p-3}$ and
$A_{2p+2}=(\lambda-2)A_{2p}-A_{2p-2}$.

$\textbf{D})$  $\mathcal{D}_{n}(\lambda),$  the determinant
associated to the monoid $D^{\infty}_{n},$ has the form
$\mathcal{D}_{n}(\lambda)=\lambda^{[\frac{n-1}{2}]}D_{n}(\lambda)$
where $D_{n}$ is a polynomial of degree $[\frac{n+2}{2}];$ the
sequence $(D_{n})_{n\geq3}$ is defined by:

$1)$ $D_{3}(\lambda)=\lambda^2-3\lambda+1,
D_{4}(\lambda)=\lambda^3-4\lambda^2+3\lambda-1,$ and
$D_{5}(\lambda)=\lambda^3-5\lambda^2+6\lambda-2$.

$2)$ $D_{2p+1}=(\lambda-2)D_{2p-1}-D_{2p-3}$ and
$D_{2p+2}=(\lambda-2)D_{2p}-D_{2p-2}$.
\end{rk}

An explicit formula for $K_{*}(\lambda)$ is given by
\begin{lem}\label{Tsigma} Denote by $q(\lambda)$ the quadratic polynomial
$\lambda^{2}-4\lambda$. Then we have

 $1)$ $K_{2p}(\lambda)=\frac{1}{2^{p-1}}\sum\limits_{i=0}^{[\frac{p}{2}]}
{\small\begin{pmatrix}p\\2i\\\end{pmatrix}}(\lambda-2)^{p-2i}q^i(\lambda);$

$2)$
$K_{2p+1}(\lambda)=\frac{1}{2^{p}}\sum\limits_{i=0}^{[\frac{p}{2}]}
{\small\begin{pmatrix}p\\2i\\\end{pmatrix}}(\lambda-2)^{p-2i}q^i(\lambda)+
\frac{1}{2^{p}}(\lambda-4)\sum\limits_{i=0}^{[\frac{p-1}{2}]}
{\small\begin{pmatrix}p\\2i+1\\
\end{pmatrix}}(\lambda-2)^{p-2i-1}q^i(\lambda)$.
\end{lem}

\proof The solutions of the characteristic equation
$\alpha^{2}-(\lambda-2)\alpha+1=0$ are\\
$\alpha_{1,2}=\frac{1}{2}(\lambda-2\pm\sqrt{q(\lambda)})$. Hence
$K_{k}=h_{1}\alpha^p_{1}+h_{2}\alpha^p_{2}$ with
\begin{eqnarray*}
(h_{1},h_{2}) &=&(1,1)
\,\,\mbox{for}\,\, k=2p, \\
(h_{1},h_{2}) &=&
\Bigg(\frac{\sqrt{q(\lambda)}+\lambda-4}{2\sqrt{q(\lambda)}},\frac{\sqrt{q(\lambda)}-
\lambda+4}{2\sqrt{q(\lambda)}}\Bigg)\,\mbox{for}\,\, k=2p+1.
\end{eqnarray*}\endproof

\section* {\bf \S \,3. Chebyshev recurrence}
The recurrences of proposition \ref{SnTn} and those of remark \ref{rk2} are close to the recurrence
for classical Chebyshev polynomials,
$T_{n}(\lambda)=\cos(\arccos\lambda)$: $T_{n}(\lambda)=\lambda
T_{n-1}(\lambda)-T_{n-2}(\lambda)$ with initial values $T_{0}=1,
T_{1}=\lambda$ (see \cite{borw}). We prove, under general
hypothesis, that all the roots of such polynomials are real and
contained in a bounded common interval $[a,b]$; in the classical
case this interval is $[-1,1]$, in our case the interval is $[0,4]$.

\begin{prop}\label{Rn}
Let $a<b$ two real numbers. If the sequence of polynomials
$\{R_{n}\}_{n\geq0}$, $R_{n}\in \mathbb{R}[X]$ satisfies the
following conditions$:$

$a)$ \textbf{[degree]}: $R_{0}=$ constant, $R_{1}=$ a polynomial of
degree 1,

$b)$ \textbf{[recurrence]}: there are $\alpha,\beta\in\mathbb{R}$,
$\alpha\neq0$ such that $$R_{n}=(\alpha X+\beta)R_{n-1}-
R_{n-2},\,n\geq2,$$

$c)$ \textbf{[boundary conditions]}: for any $n\geq0$ we have either

$i)\,R_{n}(a)R_{n+1}(a)>0,R_{n}(b)R_{n+1}(b)<0$ \,or

$ii)\,R_{n}(a)R_{n+1}(a)<0,R_{n}(b)R_{n+1}(b)>0$

then the following relations hold:

$1)$ $\deg R_{n}=n,$

$2)$ all the roots of $R_{n}$ are real, distinct and contained in
$(a,b)$:
$$a< x^{(n)}_{1}< x^{(n)}_{2}<\cdots< x^{(n)}_{n}<b,$$

$3)$ the roots of $R_{n+1}$ separate the roots of $R_{n}$, i.e.,
$$a<x^{(n+1)}_{1}< x^{(n)}_{1}< x^{(n+1)}_{2}< x^{(n)}_{2}<\cdots<
x^{(n)}_{n}< x^{(n+1)}_{n+1}<b.$$
\end{prop}

\proof We prove the proposition under the hypothesis $c(i)$.
Changing the sequence $\{R_{n}\}$ with $\{-R_{n}\}$ we can suppose
that $R_{0}(a)>0$. Starting induction by $n$, let
$\{x^{(m)}_{i}\}_{i=1,\ldots,m}$ be the roots of $R_{m}$. Suppose
the roots $\{x^{(n+1)}_{i}\}$ separate the roots $\{x^{(n)}_{i}\}$
and the signs of $R_{n}\big(x^{(n+1)}_{i}\big)$ and
$R_{n+1}\big(x^{(n)}_{i}\big)$ are those given in the following
table:
\begin{center}
 \begin{picture}(320,100)
 \put(35,20){\line(1,0){280}}\put(35,50){\line(1,0){280}}\put(35,80){\line(1,0){280}}
 \put(35,17){\line(0,1){6}}\put(35,47){\line(0,1){6}}\put(35,77){\line(0,1){6}}
 \put(315,17){\line(0,1){6}}\put(315,47){\line(0,1){6}}\put(315,77){\line(0,1){6}}

\put(87,78){$\line(0,1){4}$}\put(122,78){$\line(0,1){4}$}
\put(227,78){$\line(0,1){4}$}\put(262,78){$\line(0,1){4}$}

\put(68,48){$\line(0,1){4}$}\put(104,48){$\line(0,1){4}$}\put(141,48){$\line(0,1){4}$}
\put(208,48){$\line(0,1){4}$}\put(245,48){$\line(0,1){4}$}\put(281,48){$\line(0,1){4}$}

\put(51,18){$\line(0,1){4}$}\put(87,18){$\line(0,1){4}$}\put(122,18){$\line(0,1){4}$}\put(157,18){$\line(0,1){4}$}
\put(190,18){$\line(0,1){4}$}\put(227,18){$\line(0,1){4}$}\put(262,18){$\line(0,1){4}$}\put(298,18){$\line(0,1){4}$}

\put(102,77){$\bullet$}\put(243,77){$\bullet$}

\put(85,47){$\bullet$}\put(120,47){$\bullet$}\put(225,47){$\bullet$}\put(260,47){$\bullet$}

\put(66,17){$\bullet$}\put(102,17){$\bullet$}\put(139,17){$\bullet$}\put(206,17){$\bullet$}
\put(243,17){$\bullet$}\put(279,17){$\bullet$}

\put(85,86){{\tiny$x^{(n)}_{1}$}}\put(120,86){{\tiny$x^{(n)}_{2}$}}
\put(225,86){{\tiny$x^{(n)}_{n-1}$}}\put(260,86){{\tiny$x^{(n)}_{n}$}}

\put(65,56){{\tiny$x^{(n+1)}_{1}$}}\put(102,56){{\tiny$x^{(n+1)}_{2}$}}\put(138,56){{\tiny$x^{(n+1)}_{3}$}}
\put(203,56){{\tiny$x^{(n+1)}_{n-1}$}}\put(243,56){{\tiny$x^{(n+1)}_{n}$}}\put(277,56){{\tiny$x^{(n+1)}_{n+1}$}}

\put(48,26){{\tiny$x^{(n+2)}_{1}$}}\put(85,26){{\tiny$x^{(n+2)}_{2}$}}\put(120,26){{\tiny$x^{(n+2)}_{3}$}}
\put(151,26){{\tiny$x^{(n+2)}_{4}$}}\put(187,26){{\tiny$x^{(n+2)}_{n-1}$}}\put(225,26){{\tiny$x^{(n+2)}_{n}$}}
\put(257,26){{\tiny$x^{(n+2)}_{n+1}$}}\put(288,26){{\tiny$x^{(n+2)}_{n+2}$}}

\put(86,71){{\tiny$0$}}\put(120.5,71){{\tiny$0$}}\put(225.8,71){{\tiny$0$}}\put(260.8,71){{\tiny$0$}}

\put(32.5,70){{\small+}}\put(101,72){{\small$-$}}\put(242,72){{\tiny$\mp$}}\put(312,70){{\tiny$\pm$}}

\put(32.5,40){{\small+}}\put(84,42){{\small$-$}}\put(119.5,41){{\small+}}\put(224,42){{\tiny$\mp$}}
\put(260,42){{\tiny$\pm$}}\put(312,42){{\tiny$\mp$}}

\put(32.5,10){{\small+}}\put(312,10){{\tiny$\pm$}}
\put(65,12){{\small$-$}}\put(101.5,11){{\small+}}\put(138,12){{\small$-$}}\put(205,12){{\tiny$\mp$}}
\put(243,11){{\tiny$\pm$}}\put(278,12){{\tiny$\mp$}}

\put(67,41){{\tiny$0$}}\put(103,41){{\tiny$0$}}\put(140,41){{\tiny$0$}}\put(207,41){{\tiny$0$}}
\put(243.5,41){{\tiny$0$}}\put(279.5,41){{\tiny$0$}}

\put(49.5,11){{\tiny$0$}}\put(86,11){{\tiny$0$}}\put(120.5,11){{\tiny$0$}}\put(156,11){{\tiny$0$}}
\put(189,11){{\tiny$0$}}\put(225.5,11){{\tiny$0$}}\put(260.7,11){{\tiny$0$}}\put(297,11){{\tiny$0$}}

\put(33,85){$a$}\put(315,85){$b$}
\put(33,55){$a$}\put(315,55){$b$}\put(33,25){$a$}\put(315,25){$b$}
\put(168,82){$\cdots$}\put(168,52){$\cdots$}\put(168,22){$\cdots$}
\put(5,77){$R_{n}$}\put(5,47){$R_{n+1}$}\put(5,17){$R_{n+2}$}
 \end{picture}
 \end{center}
This table is true for $n=0$, by hypothesis:

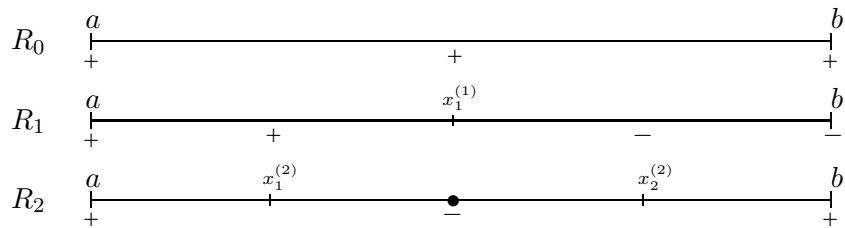
\begin{center}
 \begin{picture}(320,100)
 \put(35,20){\line(1,0){280}}\put(35,50){\line(1,0){280}}\put(35,80){\line(1,0){280}}
 \put(35,17){\line(0,1){6}}
 \put(35,47){\line(0,1){6}}
 \put(35,77){\line(0,1){6}}
 \put(315,17){\line(0,1){6}}
 \put(315,47){\line(0,1){6}}
 \put(315,77){\line(0,1){6}}
\put(101.5,42){{\small+}}\put(240,42){{\small$-$}}
\put(172,48){\line(0,1){4}} \put(168,56){{\tiny$x^{(1)}_{1}$}}
\put(100,26){{\tiny$x^{(2)}_{1}$}}\put(242,26){{\tiny$x^{(2)}_{2}$}}
\put(32.5,70){{\small+}}\put(170,72){{\small+}}\put(312,70){{\small+}}
\put(32.5,40){{\small+}} \put(312,42){{\small$-$}}
\put(32.5,10){{\small+}}\put(312,10){{\small+}}
\put(168,12){{\small$-$}} \put(169.5,17){$\bullet$}
\put(103,18){\line(0,1){4}}\put(244,18){\line(0,1){4}}
 \put(33,85){$a$}\put(315,85){$b$}
\put(33,55){$a$}\put(315,55){$b$}\put(33,25){$a$}\put(315,25){$b$}
\put(5,77){$R_{0}$}\put(5,47){$R_{1}$}\put(5,17){$R_{2}$}
 \end{picture}
 \end{center}

Using the given recurrence we get
$$R_{n+2}\big(x^{(n+1)}_{i}\big)=-R_{n}\big(x^{(n+1)}_{i}\big),$$ so we
have opposite signs on the lines $R_{n}$ and $R_{n+2}$. Now
$R_{n+2}$ has at least one root in the intervals:
$x^{(n+2)}_{1}\in\big(a,x^{(n+1)}_{1}\big),$
$x^{(n+2)}_{j}\in\big(x^{(n+1)}_{j-1},x^{(n+1)}_{j}\big),2\leq j\leq
n+1$ and $x^{(n+2)}_{n+2}\in\big(x^{(n+1)}_{n+1},b\big)$. The polynomial
$R_{n+2}$ has degree $n+2,$ hence all the  roots of
$R_{n+2}$ are $x^{(n+2)}_{1},\ldots,x^{(n+2)}_{n+2}$ and the pattern
of signs of $R_{n+2}\big(x^{(n+1)}_{i}\big)$ and
$R_{n+1}\big(x^{(n+2)}_{i}\big)$ is the same. The version $c(ii)$ is similar.
\endproof
The boundary conditions are satisfied by the polynomials $K_{n}$:

\begin{lem}\label{T(0)}
For $p\geq0:$

$1)$ $K_{2p}(0)=2(-1)^{p},\, K_{2p+1}(0)=(-1)^{p}(2p+1)\,;$

$2)$ $K_{2p}(4) =2,K_{2p+1}(4) = 1.$
\end{lem}
\proof This can be checked using recurrence or lemma \ref{Tsigma}.
\endproof

Using lemma  \ref{T(0)}, in proposition \ref{SnTn} the precise
multiplicity of the zero root of $\mathcal{K}_{n}$ is
$[\frac{n+1}{2}]$.
\begin{rk}
The same is true for $A_{n}$ but not for $D_{n}$. One can prove that
$D_{n}$ has at most two roots in $\mathbb{C}\smallsetminus
\mathbb{R};$ for example, $D_{4},$ $D_{6},$ $D_{7}$ have two complex
roots, $D_{5}$ have has only real roots.
\end{rk}

\begin{rk}
We cannot expect to obtain Chebyshev type recurrence for the
polynomials corresponding to Artin spherical monoids. For
instance, in the case the monoid $A^{+}_{n}$, Deligne theorem
[\cite{p.del}, 4.14] shows that the degree of the denominator is
{\tiny$\left(
 \begin{array}{c}
   n+1 \\
    2 \\
    \end{array}
    \right)$}
= the number of hyperplanes of the braid arrangement.
\end{rk}

Applying proposition \ref{Rn} to the sequences $(K_{2n})$ and
$(K_{2n+1})$ we obtain:
\begin{thm}
$a)$ All the roots of \,$K_{n}$ are real and belong to $(0,4)$.\\
$b)$ All the roots of $\mathcal{K}_{n}$ are real and belong to
$[0,4)$, with zero root of multiplicity $[\frac{n+1}{2}]$ and other
$[\frac{n}{2}]$ nonzero simple roots.
\end{thm}

\begin{cor}
The growth rate of $c^{[n]}_{k}$ is less than $4$ for any $n$.
\end{cor}

\proof Let $m=[\frac{n}{2}]$ and\,
$0<\alpha_{1}<\alpha_{2}<\cdots<\alpha_{m}=\alpha<4$ are the $m$
roots of $K_{n}$. From propositions \ref{HtnSn} and \ref{SnTn} we
have
\begin{eqnarray*}
 \mathcal{H}^{[n]}_{K}(t)&=&\sum\limits_{k=0}^\infty
c^{[n]}_{k}t^k=\frac{1}{t^n\mathcal{K}_{n}(\frac{1}{t})}=\frac{1}{t^mK_{n}(\frac{1}{t})}\\&=&
 \frac{1}{(1-\alpha_{1}t)(1-\alpha_{2}t)\cdots(1-\alpha_{m}t)} \\
 &=&\frac{\beta_{1}}{1-\alpha_{1}t}+\frac{\beta_{2}}{1-\alpha_{2}t}+\cdots+\frac
{\beta_{m}}{1-\alpha_{m}t}=
\sum\limits_{k=0}^\infty\Big(\sum\limits_{i=1}^m\beta_{i}(n)\alpha_{i}^k\Big)t^k
\end{eqnarray*}
\endproof

\section* {\bf \S \,4. The monoid $E^{\infty}_{n}$}

The reason to replace the standard order of vertices

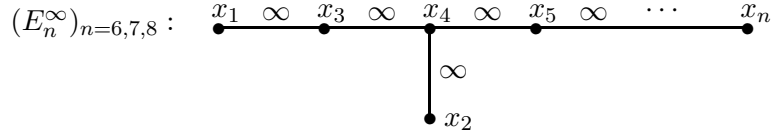
\begin{center}
\begin{picture}(310,50)
\thicklines
\put(90,37){\line(1,0){200}}\put(168.3,37){\line(0,-1){35}}
\put(86,34){$\bullet$}\put(126,34){$\bullet$}\put(166,34){$\bullet$}\put(166,0){$\bullet$}
\put(206,34){$\bullet$}\put(286,34){$\bullet$}
\put(86,41){$x_{1}$}\put(126,41){$x_{3}$}\put(166,41){$x_{4}$}\put(174,0){$x_{2}$}
\put(206,41){$x_{5}$}\put(286,41){$x_{n}$}\put(250,41){$\cdots$}
\put(105,40){$\infty$}\put(145,40){$\infty$}\put(185,40){$\infty$}
\put(225,40){$\infty$}\put(172,18){$\infty$}
\put(10,36){$(E^{\infty}_{n})_{n=6,7,8}:$}
\end{picture}
\end{center}
by the order given in the introduction is to obtain a complete
presentation (in the standard order one have to solve the ambiguity
$x_{3}x_{2}x_{1}$ and to add new relations: see \cite{berg}).

The growth rate of these monoids can be reduced to the growth rate
of $E^{\infty}_{8}$:
\begin{lem}\label{E678}
There exist homogeneous injective morphisms $E^{\infty}_{6}\rightarrowtail
E^{\infty}_{7}\!\rightarrowtail \!E^{\infty}_{8}$ defined by
$x_{i}\mapsto x_{i}$.
\end{lem}

\begin{lem}\label{Esys.rr}
In the monoid $E^{\infty}_{8}$ the following are satisfied

$a)$ $b_{0}=1$, $b_{1;i}=1;$ $b_{k} =\sum\limits_{i=1}^8
b_{k;i}\,\,\,(k\geq1);$

$c)$ $b_{k;i}\,(k\geq2)$ are given by the recurrence
 $$\left\{
  \begin{array}{ll}
    b_{k;j}=\sum\limits_{i=j-1}^8 b_{k-1;i}\,\,\,\,\,(j\neq 5),\\
    b_{k;5}=b_{k-1;3}+\sum\limits_{i=5}^8 b_{k-1;i}.
  \end{array}
\right.$$
\end{lem}
\proof (of lemma \ref{E678} and \ref{Esys.rr}):
As in lemma \ref{cpc} and proposition \ref{can.forms}, the presentation of $E^{\infty}_{*}$
is complete and the word problem in these monoids has a simple solution.
\endproof
The characteristic polynomial of this recurrence is given by:
$$\mathcal{E}_{8}(\lambda)=\lambda^4(\lambda-1)(\lambda^3-7\lambda^2+14\lambda-7)$$
having the nonzero roots in $(0,4)$:
$\lambda_{1}\approx0.75,\lambda_{2}=1,\lambda_{3}\approx2.44,\lambda_{4}\approx3.80$.

\begin{rk}
One can consider the whole sequence $(E^{\infty}_{n})_{n\geq6}:$ the
corresponding polynomials are given by
$\mathcal{E}_{n}(\lambda)=\lambda^{[\frac{n}{2}]}(\lambda-1)K_{n-1}(\lambda)$
and they satisfy the hypothesis of proposition \ref{Rn} (if we discard the common roots $0$ and $1$).
\end{rk}

\section* {\bf \S \,5. Artin spherical monoids }

\proof (of Proposition \ref{p1}) The maps
$\phi:K^{\infty}_{n}\rightarrow A^{\infty}_{n}$ and
$\phi:K^{\infty}_{n}\rightarrow D^{\infty}_{n}$, $\phi(y_{i})=x_{i}$
and $\psi:\mathcal{X}^{\infty}_{n}\rightarrow \mathcal{X}^{+}_{n}$,
$\psi(x_{i})=x_{i}$ are well defined because the images of the
generators satisfy the defining relations of the domains and these
images generate the range-monoids.\endproof

Using the notations given in the introduction:
\begin{cor}\label{a<b<c}
$a^{[n]}_{k}\leq b^{[n]}_{k}\leq c^{[n]}_{k}$.
\end{cor}

\proof (of Theorem 1: connected graphs) For the monoids
$E^{\infty}_{6},E^{\infty}_{7}$ and $E^{\infty}_{8}$ we use the
corollary \ref{a<b<c} and the results of $\S\,4$ and obtain the
inequalities
$$a^{[n]}_{k}\leq b^{[n]}_{k}\leq b^{[8]}_{k}<\beta\,\gamma^{k}<\beta\,4^{k}\,\,\textrm{where}\,\,
n=6,7,8.$$ For the other monoids we put together the inequalities
obtained before:
$$a^{[n]}_{k}\leq b^{[n]}_{k}\leq c^{[n]}_{k}<\beta(n)\,\delta^{k}<\beta(n)\,4^{k}.$$
Here $\beta=$ constant, $\beta(n)=$ constant in $k$, and $\gamma,\delta$ are between
the greatest root of the corresponding polynomials and $4$.
\endproof

The last step in the proof is to  show that an Artin spherical
monoid with a non-connected Coxeter diagram has growth rate less
than 4. The monoid associated to a disjoint union of graphs
$\Gamma_{1}\sqcup\Gamma_{2}$ is the direct product of the monoids
associated to $\Gamma_{1}$ and $\Gamma_{2}$.

\begin{lem}\label{non}
If $\mathcal{M}_{1}$ and $\mathcal{M}_{2}$ are two monoids with
growth rate less than $\gamma$, then
$\mathcal{M}_{1}\times\mathcal{M}_{2}$ has growth rate less than
$\gamma$.
\end{lem}
\proof Let us denote by $m_{1}(k)$, $m_{2}(k)$, $p(k)$, the numbers
of words of length $k$ $(k\geq1)$ in  $\mathcal{M}_{1}$,
$\mathcal{M}_{2}$ and $\mathcal{M}_{1}\times\mathcal{M}_{2}$
respectively. Then
$$p(k)=m_{1}(k)+m_{1}(k-1)m_{2}(1)+\cdots+m_{2}(k).$$
Choose $\delta>1$ and $\gamma_{0}$ such that for $i=1,2,$
$$\overline{\lim\limits_k}\,\exp\Bigg(\frac{\log m_{i}(k)}{k}\Bigg)<\delta<\gamma_{0}<\gamma,$$
and two constants $c_{1}$ and $c_{2}$ such that for any
$k$, $m_{1}(k)\leq c_{1}\delta^{k}$ and $m_{2}(k)\leq
c_{2}\delta^{k}$. Take $c_{3}=\max\{c_{1},c_{2},c_{1}c_{2}\}$ and we
obtain $p(k)\leq(k+1)c_{3}\delta^{k}$. We can find another constant
$c$ such that $(k+1)\,c_{3}\,\delta^{k}<c\gamma^{k}_{0}$. Now
$p(k)<(k+1)c\gamma^{k}_{0}$, hence $p(k)$ has growth rate less than
$\gamma$.
\endproof
\proof (of Theorem \ref{th1}: non-connected graphs) Take $\gamma=4$
and apply the lemma inductively for a finite number of connected
Coxeter graphs.
\endproof
\proof (of Theorem \ref{th2})
Garside's decomposition
$\mathcal{X}_{n}=\bigcup\limits_{k\in\mathbb{Z}}\Delta_{n}^{k}\,\mathcal{X}^{+}_{n}$
gives a nonincreasing length map
$\varphi:\,<\nabla>\times\, \mathcal{X}_{n}^{+}\rightarrow\mathcal{X}_{n}$ defined by
$\varphi(\nabla^{k},m)=\Delta^{-k}_{n}m$ (in the codomain
we consider length function using monoid generators $\nabla=\Delta_{n}^{-1}$ and
$x_{1},\cdots,x_{n})$. This map is surjective because Garside element $\Delta_{n}$
has two basic properties:

- any generator $x_{i}$ is a (left) divisor of $\Delta_{n}$ in the monoid
$\mathcal{X}_{n}^{+}:\Delta_{n}=x_{i}w_{i},$ hence $x_{i}^{-1}=w_{i}\Delta_{n}^{-1}=w_{i}\nabla$;

- $\Delta_{n}$ commutes with $x_{1},\cdots,x_{n}$ up to a permutation:
$\Delta_{n}x_{i}=x_{\sigma(i)}\Delta_{n}$, therefore any element $g\in\mathcal{X}_{n}$
can be written as a product $g=\nabla^{k}m$, where $k\geq0$ and $m\in\mathcal{X}_{n}^{+}$.
Moreover, take $g\in\mathcal{X}_{n}$ of length $l$ in the (monoid) alphabet
$\{x_{1},\cdots,x_{n},\nabla\}$ : $g=\nabla^{k_{1}}m_{1}\nabla^{k_{2}}m_{2}\cdots\nabla^{k_{s}}m_{s}$,
where $l=\sum k_{i}+\sum\mbox{length}(m_{i})$; commuting with $\nabla$, $g$ can be written as the product
$g=\nabla^{\sum k_{i}}m'_{1}m'_{2}\cdots m'_{s}=\nabla^{k}m=\varphi(\Delta^{k},m)$
and the length of $g$, length of $\nabla^{k}m$ (in the monoid $\mathcal{X}_{n}$)
and length of $(\nabla^{k},m)$ (in $M(\nabla)\times\mathcal{X}_{n}^{+}$) are equal ($m'_{i}$ is
$m_{i}$ with permuted generators $x_{j}$).
Applying lemma \ref{non} we obtain the result.
\endproof

\end{document}